\documentclass[12pt]{amsart}
\usepackage{amsmath,amssymb,latexsym}
\def\Pcl{\Sigma\mathrm{cl}}
\def\dcl{\mathrm{dcl}}
\def\cl{\mathrm{cl}}
\def\bdd{\mathrm{bdd}}
\def\Cb{\mathrm{Cb}}

\def\M{\mathfrak M}

\def\tp{\mathrm{tp}}
\def\lstp{\mathrm{lstp}}
\def\Ind#1#2{#1\setbox0=\hbox{$#1x$}\kern\wd0\hbox to 0pt{\hss$#1\mid$\hss}
\lower.9\ht0\hbox to 0pt{\hss$#1\smile$\hss}\kern\wd0}
\def\ind{\mathop{\mathpalette\Ind{}}}
\def\Notind#1#2{#1\setbox0=\hbox{$#1x$}\kern\wd0\hbox to 0pt{\mathchardef
\nn="3236\hss$#1\nn$\kern1.4\wd0\hss}\hbox to 0pt{\hss$#1\mid$\hss}\lower.9\ht0
\hbox to 0pt{\hss$#1\smile$\hss}\kern\wd0}

\newtheorem{definition}{Definition}
\newtheorem{theorem}{Theorem}
\newtheorem{proposition}[theorem]{Proposition}
\newtheorem{fact}[theorem]{Fact}
\newtheorem{lemma}[theorem]{Lemma}
\newtheorem{corollary}[theorem]{Corollary}
\newtheorem{remark}[theorem]{Remark}
\newtheorem*{claim}{Claim}
\newtheorem*{expl}{Example}
\def\bsp{\begin{expl}}
\def\ebsp{\end{expl}}
\def\beh{\begin{claim}}
\def\ebeh{\end{claim}}
\def\defn{\begin{definition}\upshape}
\def\edefn{\end{definition}}
\def\satz{\begin{theorem}}
\def\esatz{\end{theorem}}
\def\tats{\begin{fact}}
\def\etats{\end{fact}}
\def\kor{\begin{corollary}}
\def\ekor{\end{corollary}}
\def\lmm{\begin{lemma}}
\def\elmm{\end{lemma}}
\def\bem{\begin{remark}}
\def\ebem{\end{remark}}
\def\bew{\par\noindent{\em Proof: }}

\def\satzli{\begin{proposition}}
\def\esatzli{\end{proposition}}

\begin{document}
\title{Some remarks on one-basedness}
\author{Frank O. Wagner}
\address{Frank O~Wagner, Institut Girard Desargues, Universit\'e
Claude Bernard (Lyon-1), Math\'ematiques, 21 avenue Claude Bernard,
69622 Villeurbanne-cedex, France}
\email{wagner@desargues.univ-lyon1.fr}
\date{21 January 2003}
\thanks{I should like to thank Zo\'e Chatzidakis for fruitful discussions}
\begin{abstract}A type analysable in one-based types in a simple theory
is itself one-based.\end{abstract}
\maketitle

\section{Introduction}
Recall that a type $p$ over a set $A$ in a simple theory is {\em
one-based} if for any tuple $\bar a$ of realizations of $p$ and any
$B\supseteq A$ the canonical base $\Cb(\bar a/B)$ is contained in
$\bdd(\bar aA)$. One-basedness implies that the forking geometry is
particularly well-behaved; for instance one-based groups are
bounded-by-abelian-by-bounded. Ehud Hrushovski showed in
\cite[Proposition 3.4.1]{udi} that for stable stably embedded types
one-basedness is preserved under analyses: If $p$ is stable stably
embedded in a supersimple theory, and analysable (in the technical
sense defined in the next section) in one-based types, then $p$ is
itself one-based. Zo\'e Chatzidakis then gave another 
proof for supersimple structures \cite[Theorem 3.10]{zoe}, using
semi-regular analyses. We shall give an easy direct proof of the theorem
stated in the abstract, thus removing the hypotheses of stability,
stable embedding, or supersimplicity; it is similar to Hrushovski's
proof, but does not use germs of definable functions (which work less
well in simple unstable theories), and has to deal with
non-stationarity of types. While we are at it, we shall
also generalize the notion of bounded closure and
one-basedness to $\Sigma$-closure and $\Sigma$-basedness, where
$\Sigma$ is an $\emptyset$-invariant collection of partial types
(thought of as small). This may for instance be applied to consider
one-basedness {\em modulo types of finite $SU$-rank}, or {\em modulo
superstable types}.

Our notation is standard and follows \cite{wa00}.
Throughout the paper, the ambient theory will be simple, and we shall
be working in $\M^{heq}$, where $\M$ is a suffiviently saturated model
of the ambient theory. Thus tuples are tuples of hyperimaginaries, and
$\dcl=\dcl^{heq}$.

\section{$\Sigma$-closure}
In this section $\Sigma$ will be an $\emptyset$-invariant family of
partial types. We first recall the notions of internality and
analysability. 
\defn Let $\pi$ be a partial type over $A$. Then $\pi$ is\begin{itemize}
\item ({\em almost,} resp.)\ {\em internal\/} in $\Sigma$, or ({\em
almost,} resp.)\ {\em $\Sigma$-internal}, if for every realization $a$
of $\pi$ there is $B\ind_Aa$ and $\bar b$ realizing types in $\Sigma$
based on $B$, such that $a\in\dcl(B\bar b)$ (or $a\in\bdd(B\bar b)$,
respectively). 
\item {\em analysable\/} in $\Sigma$, or {\em $\Sigma$-analysable},
if for any $a\models\pi$ there are
$(a_i:i<\alpha)\in\dcl(A,a)$ such that $\tp(a_i/A,a_j:j<i)$ is
$\Sigma$-internal for all $i<\alpha$, and
$a\in\bdd(A,a_i:i<\alpha)$.\end{itemize} 
A type $\tp(a/A)$ is {\em foreign} to $\Sigma$ if $a\ind_{AB}\bar b$
for all $B\ind_Aa$ and $\bar b$ realizing types in $\Sigma$ over
$B$.\edefn 
\defn The {\em $\Sigma$-closure\/} $\Pcl(A)$ of a set $A$ is the
collection of all hyperimaginaries $a$ such that $\tp(a/A)$ is
$\Sigma$-analysable.\edefn 
We think of $\Sigma$ as a family of small types. For instance, if
$\Sigma$ is the family of all bounded types, then
$\Pcl(A)=\bdd(A)$. Other possible choices might be the family of all
types of $SU$-rank $<\omega^\alpha$, for some ordinal $\alpha$, or the
family of all superstable types. If $P$ is an $\emptyset$-invariant
family of types, and $\Sigma$ is the family of all $P$-analysable
types to which all types in $P$ are foreign, then $\Pcl(A)=\cl_P(A)$ as
defined in \cite[Definition 3.5.1]{wa00}; if $P$ consists of a single
regular type $p$, this in turn is the $p$-closure from \cite{udi1}
(see also \cite[p.\ 265]{pill}). 
\bem In general $\bdd(A)\subseteq\Pcl(A)$; if the
inequality is strict, then $\Pcl(A)$ has the same cardinality as the
ambient monster model, and hence violates the usual
conventions. However, this is usually harmless. Note that
$\Pcl(.)$ is a closure operator.\ebem 
\tats\label{forequ} The following are equivalent:\begin{enumerate}
\item $\tp(a/A)$ is foreign to $\Sigma$. 
\item $a\ind_A\Pcl(A)$.
\item $a\ind_A\dcl(aA)\cap\Pcl(A)$.
\item $\dcl(aA)\cap\Pcl(A)\subseteq\bdd(A)$.\end{enumerate}\etats
\bew This follows immediately from \cite[Proposition 3.4.12]{wa00}; see
also \cite[Lemma 3.5.3]{wa00}.\qed

$\Sigma$-closure is well-behaved with respect to independence.
\lmm\label{clind} Suppose $A\ind_BC$. Then $\Pcl(A)\ind_{\Pcl(B)}\Pcl(C)$.
More precisely, for any $A_0\subseteq\Pcl(A)$ we have
$A_0\ind_{B_0}\Pcl(C)$, where
$B_0=\dcl(A_0B)\cap\Pcl(B)$. In particular,
$\Pcl(AB)\cap\Pcl(BC)=\Pcl(B)$.\elmm
\bew Let $B_1=\Pcl(B)\cap\dcl(BC)$. Then $C\ind_BA$ implies
$C\ind_{B_1}A$, and $\tp(C/B_1)$ is foreign to $\Sigma$ by Fact
\ref{forequ}~$(3\Rightarrow1)$. Hence $C\ind_{B_1}\Pcl(A)$, and
$C\ind_{B_1}A_0$.

Since $\tp(A_0/B_0)$ is foreign to $\Sigma$ by Fact \ref{forequ}, we
obtain $A_0\ind_{B_0}\Pcl(B_0)$. But
$\Pcl(B_0)=\Pcl(B)\supseteq B_1$, whence $A_0\ind_{B_0}C$ by
transitivity, and finally $A_0\ind_{B_0}\Pcl(C)$ by foreignness to
$\Sigma$ again.\qed

\section{$\Sigma$-basedness}
Again, $\Sigma$ will be an $\emptyset$-invariant family of partial types.
\defn A type $p$ over $A$ is {\em $\Sigma$-based} if $\Cb(\bar
a/\Pcl(B))\subseteq\Pcl(\bar aA)$ for any tuple $\bar a$ of
realizations of $p$ and any $B\supseteq A$.\edefn
\bem Equivalently, $p\in S(A)$ is $\Sigma$-based if $\bar
a\ind_{\Pcl(\bar aA)\cap\Pcl(B)}\Pcl(B)$ for any tuple $\bar a$ of
realisations of $p$ and any $B\supseteq A$.\ebem
\lmm\label{clbase} Suppose $\tp(a/A)$ is $\Sigma$-based, 
$A\subseteq B$, and $a_0\in\Pcl(\bar aB)$, where $\bar a$ is a tuple
of realizations of $\tp(a/A)$. Then $\tp(a_0/B)$ is
$\Sigma$-based.\elmm
\bew Let $\bar a_0$ be a tuple of realizations of $\tp(a_0/B)$, and
$C\supseteq B$. There is a tuple $\tilde a$ of realizations of
$\tp(a/A)$ such that $\bar a_0\in\Pcl(\tilde aB)$; we may choose it
such that $\tilde a\ind_{\bar a_0B}C$. Then $\Pcl(\tilde
aB)\cap\Pcl(C)\subseteq\Pcl(\bar a_0B)$ by Lemma \ref{clind}.

Put $X=\Cb(\tilde a/\Pcl(C))$. By $\Sigma$-basedness of $\tp(a/A)$ we
have $$X\subseteq\Pcl(\tilde aA)\cap\Pcl(C)\subseteq\Pcl(\bar a_0B).$$
As $\tilde a\ind_X\Pcl(C)$ we get $\tilde aB\ind_{XB}\Pcl(C)$, and hence
$\bar a_0\ind_Y\Pcl(C)$ by Lemma \ref{clind}, where
$Y=\Pcl(XB)\cap\dcl(\bar a_0XB)$. As $Y\subseteq\Pcl(C)$, we have
$$\Cb(\bar a_0/\Pcl(C))\subseteq Y\subseteq\Pcl(XB)\subseteq\Pcl(\bar
a_0B).\qed$$
\lmm\label{unionbase} If $\tp(a)$ and $\tp(b)$ are $\Sigma$-based, so
is $\tp(ab)$.\elmm
\bew Let $\bar a$ and $\bar b$ be tuples of realizations of $\tp(a)$
and $\tp(b)$, respectively, and consider a set $A$ of parameters. We
add $\Pcl(\bar a\bar b)\cap\Cb(\bar a\bar b/\Pcl(A))$ 
to the language. By $\Sigma$-basedness of $\tp(a)$ we get
$$\Cb(\bar a/\Pcl(A))\subseteq\Pcl(\bar a)\cap\Cb(\bar a\bar
b/\Pcl(A))=\dcl(\emptyset),$$
whence $\bar a\ind\Pcl(A)$; similarly $\bar b\ind\Pcl(A)$. 

Put $b_1=\Cb(\bar b/\Pcl(\bar aA))$, and choose $\bar
a'A'\models\tp(\bar aA/b_1)$ with $\bar a'A'\ind_{b_1}\bar a\bar
bA$. Then $b_1\in\Pcl(\bar a'A')$; by $\Sigma$-basedness of $\tp(a)$
and Lemma \ref{clbase} applied to $\bar ab_1\in\Pcl(\bar a\bar a'A')$
we have $\Cb(\bar ab_1/\Pcl(AA'))\subseteq\Pcl(\bar ab_1A')$.

If $Y=\Pcl(\emptyset)\cap\dcl(b_1)$, then $A\ind_Yb_1$ by Lemma
\ref{clind}, as $b_1\in\Pcl(\bar b)$ by $\Sigma$-basedness of
$\tp(b)$ and because $\bar b\ind\Pcl(A)$; since $\tp(A'/b_1)=\tp(A/b_1)$ we
also have $A'\ind_Yb_1$,
whence $A'\ind_Y\bar ab_1A$, and $A'\ind_{YA}\bar ab_1$. As
$\Pcl(YA)=\Pcl(A)$, Lemma \ref{clind} implies
$$\begin{aligned}\Cb(\bar ab_1/\Pcl(A))&=\Cb(\bar ab_1/\Pcl(AA'))
\subseteq\Pcl(\bar ab_1A')\cap\Pcl(A)\\
&\subseteq\Pcl(\bar ab_1)\subseteq\Pcl(\bar a\bar b),\end{aligned}$$
by Lemma \ref{clind} since $A'\ind_{\bar ab_1Y}A$. On the other hand,
put $C=\Cb(\bar ab_1/\Pcl(A))$. Then $\bar b\ind_{b_1}\Pcl(\bar aA)$
by definition of $b_1$, whence $\bar a\bar b\ind_{\bar
ab_1}\Pcl(A)$; as $\bar ab_1\ind_C\Pcl(A)$ we get $\bar a\bar
b\ind_C\Pcl(A)$, whence $\Cb(\bar a\bar b/\Pcl(A))\subseteq C$. So
$$\begin{aligned}\Cb(\bar a\bar b/\Pcl(A))&=\Cb(\bar
ab_1/\Pcl(A))\cap\Cb(\bar a\bar b/\Pcl(A))\\
&\subseteq\Pcl(\bar a\bar b)\cap\Cb(\bar a\bar
b/\Pcl(A))=\dcl(\emptyset),\end{aligned}$$
whence $\bar a\bar b\ind\Pcl(A)$.\qed
\kor\label{limit} If $\tp(a_i)$ is $\Sigma$-based for all $i<\alpha$, so is
$\tp(\bigcup_{i<\alpha}a_i)$.\ekor 
\bew We use induction on $\beta$ to show that
$\tp(\bigcup_{i<\beta}a_i)$ is $\Sigma$-based, for $\beta\le\alpha$. This
is clear for $\beta=0$; it follows from Lemma \ref{unionbase} for
successor ordinals. And if $\beta$ is a limit ordinal, then for any
set $A$
$$\Cb(\bigcup_{i<\beta}a_i/\Pcl(A))=\bigcup_{i<\beta}\Cb(\bigcup_{j\le
i}a_i/\Pcl(A))\subseteq\Pcl(\bigcup_{i<\beta}a_i).\qed$$
\lmm\label{indbase} If $\tp(a/A)$ is $\Sigma$-based and $a\ind A$,
then $\tp(a)$ is $\Sigma$-based.\elmm
\bew Let $\bar a$ be a tuple of realizations of $\tp(a)$, and consider
a set $B$ of parameters. For every $a_i\in\bar a$ choose $A_i$ with
$\tp(a_iA_i)=\tp(aA)$ and $A_i\ind_{a_i}(\bar a,B,A_j:j<i)$. As
$A_i\ind a_i$ we obtain $A_i\ind(\bar a,B,A_j:j<i)$, whence
$A_i\ind_{(A_j:j<i)}\bar aB$, and inductively $(A_j:j\le i)\ind\bar
aB$. Put $\bar A=\bigcup_{a_i\in\bar a}A_i$; we just saw that $\bar A\ind\bar
aB$. Now $\tp(a_i/\bar A)$ is $\Sigma$-based for all $a_i\in\bar a$,
and so is $\tp(\bar a/\bar A)$ by Corollary \ref{limit}. As $\bar
a\ind_B\bar A$, Lemma \ref{clind} implies
$$\Cb(\bar a/\Pcl(B))=\Cb(\bar a/\Pcl(\bar AB))\subseteq\Pcl(\bar
a\bar A)\cap\Pcl(B)=\Pcl(\bar a)\cap\Pcl(B),$$ 
where the last equality follows from $\bar aA\ind_{\bar a}B$ and Lemma
\ref{clind}.\qed
\kor\label{internal} If $p$ is almost internal in $\Sigma$-based types,
then $p$ is $\Sigma$-based.\ekor 
\bew Suppose $p=\tp(a/A)$, and choose $B\ind_Aa$ and $\bar b$ such
that $a\in\bdd(B\bar b)$ and $\tp(b/B)$ is $\Sigma$-based for all $b\in\bar
b$. Then $\tp(\bar b/AB)$ is $\Sigma$-based by Lemma \ref{limit}, as
is $\tp(a/AB)$ by Lemma \ref{clbase}, and $\tp(a/A)$ by Lemma
\ref{indbase}.\qed
\lmm\label{succ} If $\tp(a)$ and $\tp(b/a)$ are $\Sigma$-based, so is
$\tp(ab)$.\elmm
\bew Consider a tuple $\bar a\bar b$ of realizations of $\tp(ab)$, and
a set $A$ of parameters. As $\tp(\bar a)$ and $\tp(\bar b/\bar a)$ are
both $\Sigma$-based, we may suppose $a=\bar a$ and $b=\bar b$. Put
$C=\Cb(ab/\Pcl(A))$; again we add $\Pcl(ab)\cap C$ to the language. By
$\Sigma$-basedness of $\tp(a)$ we get $a\ind\Pcl(A)$.

Consider a Morley sequence $(a_ib_i:i<\omega)$ in $\lstp(ab/C)$; we
may assume that $(a_ib_i:i<\omega)\ind_CabA$. Since
$(a_i:i<\omega)\ind C$ we get $ab\ind(a_i:i<\omega)$. Moreover, as
$\tp(ab/C)$ is foreign to $\Sigma$, we have
$ab\ind_C\Pcl(a_ib_i:i<\omega)$. On the other hand
$C\in\dcl(a_ib_i:i<\omega)$, whence
$$C=\Cb(ab/\Pcl(a_ib_i:i<\omega)).$$
Put $b'=\Cb(ab/\Pcl(a,a_ib_i:i<\omega))$.
Then $a\in b'$, and $b'\in\Pcl(ab)$ by $\Sigma$-basedness of $\tp(b/a)$. Put
$X=\Pcl(\emptyset)\cap\dcl(b')$. Then $b'\ind_X(a_i:i<\omega)$ by
Lemma \ref{clind}; as $\tp(b'/a_i:i<\omega)$ is
$\Sigma$-based by Lemma \ref{clbase} and Corollary \ref{limit} applied
to $b'\in\Pcl(a,a_ib_i:i<\omega)$, so is $\tp(b'/X)$ by Lemma
\ref{indbase}. Put $C'=\Cb(b'/\Pcl(a_ib_i:i<\omega))$, then
$C'\subseteq\Pcl(b')\subseteq\Pcl(ab)$ by $\Sigma$-basedness.

Now $ab\ind_{b'}\Pcl(a_ib_i:i<\omega)$ by definition of $b'$; as
$b'\ind_{C'}\Pcl(a_ib_i:i<\omega)$ by definition, we get
$ab\ind_{C'}\Pcl(a_ib_i:i<\omega)$, whence $C\subseteq C'$. We obtain
$$C=C'\cap C\subseteq\Pcl(ab)\cap C=\dcl(\emptyset),$$
whence $ab\ind\Pcl(A)$.\qed
\satz Let $p$ be analysable in $\Sigma$-based types. Then $p$ is
$\Sigma$-based.\esatz 
\bew Suppose $p=\tp(a/A)$. Then there is a sequence
$(a_i:i<\alpha)\subseteq\dcl(aA)$ such that $a\in\bdd(A,a_i:i<\alpha)$
and $\tp(a_i/A,a_j:j<i)$ is internal in $\Sigma$-based types for all
$i<\alpha$. So $\tp(a_i/A,a_j:j<i)$ is $\Sigma$-based for all $i<\alpha$ by
Corollary \ref{internal}; we use induction on $i$ to show that
$\tp(a_j:j<i/A)$ is $\Sigma$-based. This is clear for $i=0$ and $i=1$; by
Lemma \ref{limit} it is true for limit ordinals, and by Lemma
\ref{succ} it holds for successor ordinals.\qed
\kor If $p$ is analysable in one-based types, then $p$ is itself
one-based.\qed\ekor

\end{document}